\documentclass[12 pt]{amsart}
\usepackage{mathpazo} 
\usepackage{xcolor}
\usepackage{tensor}
\usepackage{hyperref}

\usepackage{setspace}
\usepackage{geometry}
\usepackage{amssymb}
\usepackage{graphicx}
\usepackage[mathscr]{eucal}

\usepackage[colorinlistoftodos]{todonotes}

\usepackage{enumerate}
\newenvironment{enum}{  
\begin{enumerate}[\upshape(\arabic{section}.\arabic{equation}a)] }  
{  \end{enumerate}   }
\newcommand{\itemref}[2] {{\upshape(\ref{#1}\ref{#2})}}

\theoremstyle{plain}
\newtheorem{Thm}{Theorem}
\newtheorem{Cor}[Thm]{Corollary}
\newtheorem{Prop}[Thm]{Proposition}

\newtheorem{Remark}[Thm]{Remark}
\newtheorem{Q}[Thm]{Question}
\theoremstyle{definition}

\newcommand{\qf}{\mathscr Q}
\newcommand{\belt}{\mathscr B}
\newcommand{\dual}{\mathscr D}

\newcommand{\lma}{\langle\langle}
\newcommand{\rma}{\rangle\rangle}

\newcommand{\hs}{\mathscr H}
\newcommand{\cau}{\mathbf C}
\newcommand{\ler}{\mathbf L}
\newcommand{\lerc}{\tilde{\mathbf L}}
\newcommand{\sr}{\varphi}
\newcommand{\sff}{II}
\newcommand{\hyp}{\mathfrak h}
\newcommand{\inc}{\mathscr I}
\newcommand{\gr}{\Gamma}
\newcommand{\muf}{\mu_{\sf{Fef}}}
\DeclareMathOperator{\pv}{P.V.}

\newcommand{\levi}{\mathscr L}

\renewcommand{\Re}{\operatorname{Re}}
\renewcommand{\Im}{\operatorname{Im}}
\renewcommand{\bar}{\overline}
\renewcommand{\tilde}{\widetilde}
\newcommand{\intl}{\int\limits}
\newcommand{\suml}{\sum\limits}
\newcommand{\C}{\mathbb C}
\newcommand{\R}{\mathbb R}
\newcommand{\CP}{\mathbb{CP}}
\newcommand{\Z}{\mathbb Z}
\newcommand{\dee}{\partial}
\newcommand{\deebar}{\overline\partial}
\newcommand{\st}{\,:\,}

\newcommand{\oh}{\mathscr O}

\newcommand{\lam}{\lambda}
\newcommand{\eps}{\varepsilon}
\newcommand{\eqdef}{\overset{ \text{def} }{=}}
\newcommand{\mapdef}[4]{ #1 &\to #2 \\ #3 &\mapsto #4 }
\newcommand{\w}{\wedge}
\newcommand{\inv}{^{-1}}
\numberwithin{equation}{section}

\begin{document}
\title[ Holomorphic projection and duality]
{Holomorphic projection and duality\\ for domains in complex projective space}
\author[David E. Barrett]{David E. Barrett}
\thanks{Supported in part by NSF grant number DMS-0901205.}
\address{Dept. of Mathematics\\University of Michigan
\\Ann Arbor, MI  48109-1043  USA }
\date{\today}
\email{barrett@umich.edu}
\begin{abstract}  We show that the efficiency of a natural pairing between certain projectively invariant Hardy spaces on dual  strongly $\C$-convex real hypersurfaces in complex projective space is measured by the norm of the corresponding Leray transform. 
\end{abstract}


\maketitle

\section{Introduction} 

Let $S$ be a smooth compact real hypersurface in complex projective space $\CP^n$.  $S$ is said to be {\em strongly $\C$-convex} (or {\em strongly $\C$-linearly convex}) if all complex tangent hyperplanes to $S$ lie to one side of $S$ with minimal order of contact (see \S\ref{S:sclc} below).  When $S$ is 
strongly $\C$-convex, the set of all complex tangent hyperplanes to $S$ form a smooth strongly $\C$-convex real hypersurface $S^*$ in the dual projective space $\CP^{n*}$.

There is a natural $\C$-bilinear pairing between the space of square-integrable sections over $S$ of the $n^{\text{th}}$ power of the tautological line bundle and the corresponding space of sections over $S^*$.  The natural generalization to this setting of the one-dimensional Cauchy transform is the Leray transform defining a projection operator from the $L^2$ section spaces just described onto the corresponding Hardy spaces of holomorphic boundary values.  The Leray transform is in a suitable sense self-adjoint with respect to the pairing mentioned above.

In this paper we show that the norm of the Leray transform measures the effectiveness of the induced pairing on Hardy spaces.
The sharp form of this result requires the use of specific $L^2$ norms.  The norms we use are  M\"obius-invariant (that is, invariant under automorphisms of $\CP^n$) and determined by second-order data for $S$; the norms are in fact determined by these conditions along with the duality properties (see Remark \ref{R:uniq}).

 The paper is organized as follows. In \S \ref{S:ojk} we adapt a standard treatment of line bundles on $\CP^n$ to fit the needs of the current paper.  \S \ref{S:do} contains a brief account of the one-dimensional versions (where available) of the constructions arising later in paper.  In \S \ref{S:bt} we provide statements of the main results of this paper using affine coordinates.  The M\"obius-invariant geometry of hypersurfaces is studied in \S \ref
{S:igh}.  (This section also contains supplementary material specific to the two-dimensional situation.)  \S 
\ref{S:dual} sets forth the key facts about dual hypersurfaces.  The $L^2$ norms and the $\C$-linear pairing between $L^2$ spaces on dual hypersurfaces are defined in \S\S \ref{S:bil} and \ref{S:hard}. (The pairing is connected with the Fantappi\`e transform [Fan] and coincides -- up to a constant -- with the variant given in \S 3.2 of [APS].) Finally, \S \ref{S:ler} brings the Leray transform into the picture and includes the proofs of the main results.

The theory we construct is in some respects simpler and in some respects richer in dimension two than in higher dimensions.  The paper is written  to make it reasonably easy for the reader to focus primarily on the two dimensional case,  and some readers will want to exercise this option on a first pass.

\section{\texorpdfstring
{The bundles $\oh(j,k)$}{O(j,k)}} \label{S:ojk}
(Compare [GrHa, \S 1.3] or [APS, \S 3.2].)

We define a family $\oh(j,k)$ of $\C$-bundles  over complex projective space $\CP^n$.  (The bundles will be holomorphic only when $k=0$.)

A section of $\oh(j,k)$ over a subset $E$ of $\CP^n$
is given by a complex-valued function $F$ on the corresponding dilation-invariant subset of $\C^{n+1}\setminus\{0\}$ satisfying the homogeneity condition
\begin{equation*}
F(\lam \zeta)=\lam^j \bar{\lam^k} F(\zeta).
\end{equation*}
(In this paper, $j$ and $k$ will generally be integers but in fact it suffices to have $j,k\in\R$, $j-k\in\Z$.)

We denote by $\Gamma(E;j,k)$ the space of continuous sections of $\oh(j,k)$ over $E$.

If $F\in\Gamma(E;j,k)$ then $\bar F\in\Gamma(E;k,j)$. 
If $F_1\in\Gamma(E;j_1,k_1)$ and $F_2\in\Gamma(E;j_2,k_2)$ then $F_1F_2\in \Gamma(E;j_1+j_2,k_1+k_2)$.

It makes sense to declare that $F\in \oh(j,j)$ is {\em positive} when  $F$ takes values in $\R_+$.  If $F\in \oh(j,j)$ is positive then $F^{k/j}\in \oh(k,k)$ is well-defined and positive.

Similar remarks apply to sections which are $\ge0$.

If $F\in\Gamma(E;j,k)$ then we may define 
\begin{equation*}
|F|=\left(F\bar F\right)^{1/2}\in\Gamma\left(E;\frac{j+k}{2},\frac{j+k}{2}\right).
\end{equation*}

The bundle $\oh(-n-1,0)$ may be identified with the canonical bundle of $(n,0)$-forms by identifying  a form written as \[f(z_1,\dots,z_n) \,dz_1\w\dots\w dz_n\] in standard affine coordinates with $F\in \Gamma(E;-n-1,0)$ via the formulae
\begin{align*}
F(\zeta_0,\zeta_1,\ldots,\zeta_n)&=\zeta_0^{-n-1}f(\zeta_1/\zeta_0,\dots,\zeta_n/\zeta_0)\\
f(z_1,\dots,z_n) &= F(1,z_1,\ldots,z_n).
\end{align*}

Similarly, a volume form \[f(z_1\dots,,z_n) \,dz_1\w \dots \w dz_n\w \bar{d{z}_1\w \dots \w d{z_n}}\] may be identified with a section $F\in\Gamma(E;-n-1,-n-1)$.  For $n$ even our notion of positivity of $F$ coincides with the usual notion of positivity of the corresponding  volume form; thus $F\ge0$ implies $\int\limits_E f(z_1\dots,,z_n) \,\,dz_1\w \dots \w dz_n\w \bar{d{z}_1\w \dots \w d{z_n}}\ge0$.  For $n$ odd we have instead that $F\ge0$ implies $\int\limits_E if(z_1\dots,z_n) \,\,dz_1\w \dots \w dz_n\w \bar{d{z}_1\w \dots \w d{z_n}}\ge0$.  

For $F\in\Gamma(E;j,k)$ the above remarks allow us to interpret $\left(F\bar F\right)^{-\frac{n+1}{j+k}}=|F|^{-2\frac{n+1}{j+k}}$ as a volume form; it is guaranteed that  $\left(F\bar F\right)^{-\frac{n+1}{j+k}}\ge0$.

Automorphisms of $\CP^n$ are induced by matrices $M\in SL(n+1,\C)$.  The pullback operation $M^*:F\mapsto F\circ M$ induces lifted automorphisms of the line bundles $\oh(j,k)$ respecting the conjugation and multiplication operations.  The lifted automorphisms are in general not unique: each automorphism of $\CP^n$ is represented by by $n+1$ distinct choices of $M$ differing by roots of unity, and these give rise to $\frac{n+1}{\gcd\{j-k,n+1\}}$ distinct lifted automorphisms of $\oh(j,k)$. Note that the lifts are in fact unique for the the canonical bundle $\oh(n+1,0)$ as well as for bundles of the form $\oh(j,j)$ or $\oh(j(n+1),k(n+1))$.

\begin{Remark} \label{R:oeuc}
All the bundles $\oh(j,k)$ are in fact trivial over $\C^n$ (identified with $\{(\zeta_0:\zeta_1:\cdots:\zeta_n)\in\CP^n\st \zeta_0\ne 0\}$).  In particular, a section of $\oh(j,k)$ over $E\subset\C^n$ may be identified with a scalar function on $E$ via the formulae
\begin{align*}
F(\zeta_0,\zeta_1,\ldots,\zeta_n)&=\zeta_0^{j}\bar{\zeta_0^{k}}\,f(\zeta_1/\zeta_0,\dots,\zeta_n/\zeta_0)\\
f(z_1,\dots,z_n) &= F(1,z_1,\ldots,z_n).
\end{align*}
In this notation we have
\begin{equation*}
M^*f = 
\left(M_{0,0} + M_{0,1}z_1+\dots M_{0,n}z_n\right)^j
\cdot\bar{\left(M_{0,0} + M_{0,1}z_1+\dots M_{0,n}z_n\right)^k}\cdot
\left(f\circ \Theta_M\right)
\end{equation*}
where $M=\left(M_{j,k}\right)_{j,k=0}^n$ and 
\begin{equation*}
\Theta_M(z) = 
\left(
\frac{M_{1,0} + M_{1,1}z_1+\dots M_{1,n}z_n}{M_{0,0} + M_{0,1}z_1+\dots M_{0,n}z_n},\dots,
\frac{M_{n,0} + M_{n,1}z_1+\dots M_{n,n}z_n}{M_{0,0} + M_{0,1}z_1+\dots M_{0,n}z_n}
\right).
\end{equation*}

In view of the identification of $\oh(-n-1,0)$ with the canonical bundle, it is convenient to use the notation
\begin{equation*}
f(z_1,\dots,z_n) \left(dz_1\w \dots \w dz_n\right)^{\frac{-j}{n+1}}
\bar{
\left(dz_1\w \dots \w dz_n\right)^{\frac{-k}{n+1}}
}
\end{equation*}
for sections of $\oh(j,k)$.
\end{Remark}

\section{Dimension one} \label{S:do}

Let $\gamma$ be a smooth oriented simple closed curve in the Riemann sphere $\CP^1$, and let $\Omega_+$ and $\Omega_-$ denote the components of $\CP^1\setminus\gamma$ enclosed positively and negatively, respectively, by $\gamma$. 

Using the  conventions of Remark \ref{R:oeuc}, a section   $f$ of $\oh(-1,0)$ will be notated as $f(z)\,\sqrt{dz}$, where $z$ is the standard affine coordinate on $\C\subsetneq \CP^1$.

We have the hermitian pairing
\begin{equation}\label{E:hp1}
\langle f(z)\,\sqrt{dz}, g(z)\,\sqrt{dz}\rangle
=
\int\limits_\gamma f(z)\bar{g(z)} \,|dz|
\end{equation}
on $\Gamma(\gamma;-1,0)$.  We denote the resulting Hilbert space by $L^2(\gamma;-1,0)$.  

We denote by $\hs_+(\gamma)$ and $\hs_-(\gamma)$ the Hardy spaces consisting of boundary values in $L^2(\gamma;-1,0)$ of holomorphic sections of $\oh(-1,0)$ on $\Omega_+$ and $\Omega_-$, respectively.

We also have a $\C$-bilinear pairing
\begin{equation}\label{E:bp1}
\lma f(z)\,\sqrt{dz}, g(z)\,\sqrt{dz} \rma
=
\int\limits_\gamma f(z)g(z) \,dz
\end{equation}
on $L^2(\gamma;-1,0)$.

We define the inner and outer Cauchy transforms $\cau_\pm$ by 
\begin{equation}\label{E:CauchyTransform}
\cau_\pm \left(f(w)\,\sqrt{dw}\right) 
= \frac12 f(z)\,\sqrt{dz}\pm \frac{1}{2\pi i}\pv
\left\langle\left\langle f(w)\,\sqrt{dw}, \frac{\sqrt{dw}\sqrt{dz}}{w-z}\right\rangle\right\rangle;
\end{equation}
here the pairing is taken with respect to the $w$ variable and  
$\pv$ denotes the principal value of the singular integral.  The Cauchy transforms define bounded projection operators
\begin{equation*}
\cau_\pm: L^2(\gamma;-1,0) \to \hs_\pm(\gamma).
\end{equation*}
We have
\begin{equation}\label{E:dsa1}
\lma \cau_+ f, g \rma = \lma f, \cau_- g \rma
=\lma \cau_+ f, \cau_- g \rma.
\end{equation}
The norm $\|\cau_+\|=\|\cau_-\|$ measures the effectiveness of the pairing $\lma\cdot,\cdot\rma$ between $\hs_+(\gamma)$ and 
$\hs_-(\gamma)$:
\begin{equation}\label{E:pc1}
\frac{1}{\|\cau_\pm\|}=
\inf_{\substack{f\in \hs_+(\gamma)\\ \| f\|=1}}\;
\sup_{\substack{g\in \hs_-(\gamma)\\ \| g\|\le 1}}\;
\left|\lma f,g \rma\right|
=
\inf_{\substack{g\in \hs_-(\gamma)\\ \| g\|= 1}}\;
\sup_{\substack{f\in \hs_+(\gamma)\\ \| f\|\le 1}}\;
\left|\lma f,g \rma\right|.
\end{equation}
(The proof follows that of Corollary \ref{C:pc} below.)
The norm $\|\cau_\pm\|$ will equal 1 if and only if $\gamma$ is a circle (or extended line). (This follows from [KeSt2, \S 7].)

The constructions above are all invariant under the M\"obius group of automorphisms of $\CP^1$.  As explained in \S \ref{S:ojk}, there is a $\pm$ amibiguity in the lifting of an automorphism of $\CP^1$ to the bundle $\oh(-1,0)$, but the $\pm$ signs all cancel in formulae such as \eqref{E:hp1}, \eqref{E:bp1}, \eqref{E:dsa1} and \eqref{E:pc1}.

There are a number of basic one-dimensional results that do not admit higher-dimen\-sional versions in the theory developed below.  They include
\begin{itemize}
\item the formula $\lma f, g \rma=0$ for $f,g\in \hs_+(\gamma)$ or $f,g\in \hs_-(\gamma)$ (following from Cauchy's theorem)
\end{itemize}
and
\begin{itemize}
\item the identity $f=\cau_+ f - \cau_-f$ exhibiting  $L^2(\gamma;-1,0)$ as the algebraic direct sum of $\hs_+(\gamma)$ and $\hs_-(\gamma)$.
\end{itemize}

\section{Statement of main results (in affine coordinates)} \label{S:bt}

Here we state some of the main results of this paper in standard affine coordinates $z_1,\dots,z_n$ on $\C^n\subset\CP^n$.  (See \S\S \ref{S:bil},  \ref{S:hard}  and \ref{S:ler} and below for the projective coordinate versions.)

Let $S\subset\CP^n$ be a strongly pseudoconvex real hypersurface with defining function $\rho$.  $S$ is said to be {\em strongly $\C$-convex} if $S$ is locally M\"obius equivalent to a strongly convex domain.  (This condition is equivalent to the condition on complex tangent hyperplanes stated in the introduction -- see \S\ref{S:sclc} below.)

\begin{Thm}\label{T:aff-norm} Let $S$ and $\rho$ be as above.

\refstepcounter{equation}\label{N:an}
\begin{enum}
\item\label{I:aff-norm-def}
The formula 
\begin{equation*}\label{E:aff-norm-def}
\|F\|^2_S\eqdef
2\intl_{S\cap \C^n}
\left|
{\det} \begin{pmatrix}
0 & \rho_{\bar k}\\
\rho_j & \rho_{j,\bar k}
\end{pmatrix}\right|
\cdot\left|\det
\begin{pmatrix}
0 &0 &  \rho_k & 0\\
0 & 0  & 0 &\rho_{\bar k}\\
\rho_j & 0& \rho_{j,k} & \rho_{j,\bar k}\\
0 & \rho_{\bar j}  &\rho_{\bar j,k} &\rho_{\bar j,\bar k}\\
\end{pmatrix} 
\right|^{\frac{-n}{2(n+1)}}
\cdot
\frac{|f|^2\,dS}{\|d\rho\|}
\end{equation*}
defines a M\"obius-invariant norm on $F\in\Gamma(S;j,k)$ for $j+k=-n$.  Here $\|d\rho\|$ is the euclidean length of $d\rho$, $dS$ is euclidean surface area on $S$, subscripts denote differentiation and $f$ is the affine coefficient of $F$ as set forth in Remark
\ref{R:oeuc} above.  The norm is independent of the choice of defining function $\rho$.
\medskip
\item\label{I:no-norm} No such M\"obius-invariant norm exists for  $F\in\Gamma(S;j,k)$ when $j+k\ne-n$.
\end{enum}
\end{Thm}

See \S \ref{S:hard} for the proof of these results, including  a more detailed statement of part \itemref{N:an}{I:no-norm} and alternate descriptions of $\|F\|_S$.

The norm $\|F\|_S$ defined above is in fact not the only M\"obius-invariant norm on $\Gamma(S;j,k)$ but it is determined by the invariance property together with its duality properties (in particular, by item \itemref{N:bt}{I:bt-xd} from Theorem \ref{T:bt} below for $\Gamma(S;-n,0)$). (See Remark \ref{R:uniq}  for details.) Thus this norm is a natural object despite its slightly cumbersome appearance. 

Using $\|F\|_S$ we extend $\Gamma(S;j,k)$ to a M\"obius-invariant Hilbert space $L^2(S;j,k)$ when $j+k=-n$. When  $S$ is the boundary of a strongly pseudoconvex domain $\Omega$ we may define 
the Hardy space $\hs (S)$ to be 
$L^2(S;-n,0)$-closure of
\begin{equation} \label{E:aom}
\left\{ F\in \Gamma(\Omega\cup S;-n,0)\st F\text{ holomorphic on }\Omega\right\}.
\end{equation}  
(In this setting, this construction will agree with other standard definitions of the Hardy space.) 

When $S$ is smooth and strongly $\C$-convex then the dual hypersurface $S^*\subset\CP^{n*}$ consisting of all complex tangent hyperplanes to $S$ is also smooth and strongly $\C$-convex and satisfies $S^{**}=S$.  We have the following result.

\begin{Thm}\label{T:bt}
Let $S$ be as above.  Then there is a $\C$-bilinear pairing $\lma F,G \rma=\lma F,G \rma_{S,S^*}$ between $L^2(S;-n,0)$ and $L^2(S^*;-n,0)$ with the following properties.
\refstepcounter{equation}\label{N:bt}
\begin{enum}
\item\label{I:bt-mob}  
The pairing is M\"obius-invariant; that is, 
\[\lma M^*F, M^*G \rma_{S,S^*}=\lma F,G \rma_{M(S),M(S)^*}\]
for  $M\in SL(n+1,\C)$.
\medskip
\item\label{I:bt-iso}  Suppose that $S$ is isotopic to $\tilde S$ within a region where $F$ is holomorphic and that the corresponding homotopy between $S^*$ and $\tilde S^*$ takes place in a region where $G$ is holomorphic.  Then 
$\lma F,G \rma_{S,S^*}=\lma F,G \rma_{\tilde S,\tilde S^*}.$
\medskip
\item\label{I:bt-xd}  $\lma F,G \rma$ is an exact duality pairing between $L^2(S;-n,0)$ and $L^2(S^*;-n,0)$; that is, 
\begin{align*}
\| F \|_S &=
\max_{\substack{G\in L^2(S^*;-n,0)\\ \| G\|\le 1} }\;
\left|\lma F, G \rma\right|\\
\| G \|_{S^*}&=
\max_{\substack{F\in L^2(S;-n,0)\\ \| F\|\le 1} }\;
\left|\lma F, G \rma\right|.
\end{align*}
\medskip
\item\label{I:bt-hsd} $\lma F,G \rma$ is a duality pairing between $\hs(S)$ and $\hs(S^*)$; that is, 
\begin{align*}
\inf_{\substack{F\in \hs(S)\\ \| F\|=1}}\;
\sup_{\substack{G\in \hs(S^*)\\ \| G\|\le 1}}\;
\left|\lma F,G \rma\right|&>0\\
\inf_{\substack {G\in \hs(S^*)\\ \| G\|=1}}\;
\sup_{\substack{F\in \hs(S)\\ \| F\|\le 1}}\;
\left|\lma F,G \rma\right|&>0,
\end{align*}
and so [Bab] any bounded linear functional on either space may be uniquely represented by pairing against a function in the other space.

\medskip
\item\label{I:bt-eff}  The efficiency of the duality pairing between $\hs(S)$ and $\hs(S^*)$ is given by the norm of the Leray
transform $\ler_S$ projecting $L^2(S;-n,0)$ onto $\hs(S)$; in fact we have
\begin{align}\label{E:bt-eff}
\| \ler_S \|^{-1}
& = \inf_{\substack{F\in \hs(S)\\ \| F\|=1}}\;
\sup_{\substack{G\in \hs(S^*)\\ \| G\|\le 1}}\;
\left|\lma F,G \rma\right|\notag\\
&=\inf_{\substack {G\in \hs(S^*)\\ \| G\|=1}}\;
\sup_{\substack{F\in \hs(S)\\ \| F\|\le 1}}\;
\left|\lma F,G \rma\right| \notag\\
&=\| \ler_{S^*} \|^{-1}.\notag
\end{align}
\end{enum}
\end{Thm}

\bigskip

The pairing  is defined in \S \ref{S:bil}.  For additional details (including the definition of $\ler_S$) and proofs of the items in the theorem see \S \ref{S:ler}.

\section{Invariant geometry of hypersurfaces} \label{S:igh}

\subsection{A scalar invariant} \label{S:scinv}

\begin{Prop}\label{P:srdef}
Let $S\subset\CP^n$ be a (not necessarily compact) strongly pseudoconvex hypersurface and let $\rho$ be a defining function for $S\cap\C^n$ with $d\rho$ pointing to the pseudoconcave side of $S$.  Define real-valued functions $\delta_{\rho,1}, \delta_{\rho,2}$ on $S\cap \C^n$ by
\begin{align}
\delta_{\rho,1}&=-\det
\begin{pmatrix}
0 & \rho_{\bar k}\\
\rho_j & \rho_{j,\bar k}
\end{pmatrix}
\\
\label{E:del2-def}
\delta_{\rho,2}&=(-1)^{n+1}\det
\begin{pmatrix}
0 &0 &  \rho_k & 0\\
0 & 0  & 0 &\rho_{\bar k}\\
\rho_j & 0& \rho_{j,k} & \rho_{j,\bar k}\\
0 & \rho_{\bar j}  &\rho_{\bar j,k} &\rho_{\bar j,\bar k}\\
\end{pmatrix}. 
\end{align}
(As before, the subscripts on $\rho$ denote differentiation.) 
Then the following hold.
\refstepcounter{equation}\label{N:srdef}
\begin{enum}
\item\label{I:ineq-srdef} $\delta_{\rho,1}>0$ on $S\cap\C^n$.
\item\label{I:def-srdef}
The real-valued function
\begin{equation*}
\sr_S\eqdef \delta_{ \rho,1}^{-2}\cdot \delta_{ \rho,2}
\end{equation*}
 on $S\cap\C^n$ is independent of the choice of defining function $\rho$.
\item\label{I:mob-srdef} $\sr_S$ satisfies the transformation law $\sr_{S}=M^* \sr_{M(S)}$.
\item\label{I:ext-srdef}  $\sr_S$ extends to a smooth function on all of $S$.
\end{enum}
\end{Prop}

\begin{proof}  \itemref{N:srdef}{I:ineq-srdef} is a standard consequence of strong pseudoconvexity.

If $\rho$ is replaced by $\eta\rho$ with $\eta>0$ on $S$ then using the product rule, row operations and the vanishing of $\rho$ on $S$ we find that 
\begin{align*}
\delta_{\eta\rho,1}&=\eta^{n+1} \delta_{\rho,1}\\
\delta_{\eta\rho,2}&=\eta^{2n+2} \delta_{\rho,2},
\end{align*}
which implies 
\itemref{N:srdef}{I:def-srdef}.

To tackle the remaining claims, pick a real-valued defining function for $S$ in $\CP^n$ and multiply it by $\|\zeta\|^4$ to obtain a defining function $R$ with values in $\oh(2,2)$.
Let 
\begin{align*}
\Delta_{R,1}&= -\frac{1}{4}\det 
\begin{pmatrix}
R_{j\bar k}
\end{pmatrix}   \\
\Delta_{R,2}&= \frac{(-1)^{n+1}}{9} \det
\left(
\begin{array}{r|r}
R_{jk} & R_{j\bar k}\\
\hline
R_{\bar jk} & R_{\bar j \bar k}
\end{array}
\right). 
\end{align*}
(Here $j$ and $k$ run from $0$ to $n$.)

The entries of $\begin{pmatrix}
R_{j\bar k}
\end{pmatrix} $ take values in $\oh(1,1)$, so 
$\Delta_{R,1}$ takes values in $\oh(n+1,n+1)$, and similar considerations show that $\Delta_{R,2}$ takes values in $\oh(2n+2,2n+2)$.  It is easy to check that $\Delta_{R\circ M,1}=M^*\Delta_{R,1}$ and 
$\Delta_{R\circ M,2}=M^*\Delta_{R,2}$
for  $M\in SL(n+1,\C)$.

\itemref{N:srdef}{I:ext-srdef} and \itemref{N:srdef}{I:mob-srdef} will follow now once we establish
\begin{equation}\label{E:sr-proj-def}
\sr_S= \Delta_{ R,1}^{-2}\cdot \Delta_{R,2}.
\end{equation}
In view of the M\"obius invariance of $\Delta_{R\circ M,1}$ and $\Delta_{R\circ M,2}$ and the evident translation invariance of $\sr_S$ it will suffice to check that this holds at $(1:0:\dots:0)$.
Setting $R(\zeta_0,\zeta_1,\dots,\zeta_n)=\zeta_0^2\bar{\zeta_0^2}\, r(\zeta_1/\zeta_0,\dots,\zeta_n/\zeta_0)$, plugging this in to the definitions of $\Delta_{R,1}$ and $\Delta_{R,2}$
and then setting $\zeta_0=1$ we find after a few  row operations that the affine coefficients of $\Delta_{R,1}$ and $\Delta_{R,2}$ are $\delta_{r,1}$ and $\delta_{r,2}$.  Since $r$ is a defining function for $S\cap\C^n$, \eqref{E:sr-proj-def} now follows from 
\itemref{N:srdef}{I:def-srdef}
\end{proof}

\subsection{\texorpdfstring{Strong $\C$-linear convexity}{Strong C-linear convexity}}\label{S:sclc}

\begin{Prop}\label{P:nrm-frm}
Let $S$ be a smooth strongly pseudoconvex real hypersurface in $\CP^n$ and let $p$ be a point in $S$.  Then there is an automorphism of $\CP^n$ moving $p$ to $0\in\C^n$ so that the transformed $S$ takes the form 
\begin{equation} \label{E:exp-choose-c}
v= \suml_{j,k=1}^{n-1} \alpha_{j,k} z_j \bar z_k +\Re \left(\suml_{j,k=1}^{n-1} \beta_{j,k} z_j z_k\right)+cu^2 
+ O\left(\|(z_1,\dots,z_{n-1},u)\|^3\right)
\end{equation}
near $0$, where we have set $z_n=u+iv$.
Here the real constant $c$ may be set arbitrarily, but the terms $\suml_{j,k=1}^{n-1} \alpha_{j,k} z_j \bar z_k +\Re \left(\suml_{j,k=1}^{n-1} \beta_{j,k} z_j z_k\right)$ are determined up to a scalar multiple and a $\C$-linear change of variables in $(z_1,\dots,z_{n-1})$.
\end{Prop}

\begin{proof}
After a preliminary transformation we may assume that $p=0$ with real tangent space $T_0 S= \C^{n-1}\times\R$; thus $S$ takes the form
\begin{equation*}
v= \suml_{j,k=1}^{n-1} \alpha_{j,k} z_j \bar z_k +\Re \left(\suml_{j,k=1}^{n-1} \beta_{j,k} z_j z_k\right)
+\Im\left(\suml_{j=1}^{n-1} \gamma_j z_ju \right)
+\tilde cu^2 
+ O\left(\|(z_1,\dots,z_{n-1},u)\|^3\right).
\end{equation*}
Applying the transformation 
\begin{equation*}
z \mapsto \frac{z}{1-\sum\limits_{j=1}^{n-1} \gamma_j z_j + \left( c-\tilde c\right) z_n}
\end{equation*}
we see that \eqref{E:exp-choose-c} holds.

For the invariance properties of  $q(z_1,\dots,z_{n-1})\eqdef\suml_{j,k=1}^{n-1} \alpha_{j,k} z_j \bar z_k +\Re \left(\suml_{j,k=1}^{n-1} \beta_{j,k} z_j z_k\right)$  note that a M\"obius transformation mapping $0$ to $0$ and preserving $T_0 S= \C^{n-1}\times\R$ must take the form 
\begin{equation*}
z=(z',z_n) \mapsto \frac{(M'z'+z_n p', \kappa z_n) }{1-\sum_{j=1}^{n} \tilde\gamma_j z_j}
\end{equation*}
with $M'\in GL(n-1,\C)$, $p'\in\C^{n-1}$ and $\kappa\in\R\setminus\{0\}$.
Direct computation shows that the revised $q$ terms are given by $\kappa^{-1} q(M'z')$.
(See  [Bol2, Prop. 2] for more detail.)
\end{proof}

\begin{Remark}\label{R:iso}
With a fixed choice of $q$ we see that our choice of projective coordinates at $p$ is determined up to a 
M\"obius transformation mapping $0$ to $0$ whose derivative at $0$ takes the form
\begin{equation*}
\begin{pmatrix}
 z' \\ z_n 
\end{pmatrix}
\mapsto 
\begin{pmatrix}
 M'z' + z_n p'\\ \kappa z_n
\end{pmatrix}
\end{equation*}
with $\kappa q(z')=q(M'z')$.  Since $q$ is positive definite we must have $|\det M'|^2=\kappa^{n-1}$,
\end{Remark}

The last part of the above proof shows in fact that the restriction of the second fundamental form of $S$ to the maximal complex subspace $HS=TS\cap JTS$ defines a M\"obius-invariant 
$T\CP^n/TS$-valued quadratic form $\sff$ on $HS$. The Levi-form $\levi$ (given by \linebreak $\suml_{j,k=1}^{n-1} \alpha_{j,k} z_j \bar z_k$ in the proposition) is of course just the hermitian part  of $\sff$:
\begin{equation*}
\levi(X_1,X_2)=\frac{\sff(X_1,X_2)+\sff(JX_1,JX_2)}{2}.
\end{equation*}
We denote the anti-hermitian part of $\sff$ by $\qf$; thus  
\begin{equation*}
\qf(X_1,X_2)=\frac{\sff(X_1,X_2)-\sff(JX_1,JX_2)}{2}
\end{equation*}
and $\qf$ is given by $\Re\left(\suml_{j,k=1}^{n-1} \beta_{j,k} z_j z_k\right)$ in the proposition.

The hypersurface $S$ is said to be {\em strongly $\C$-convex} if $\sff$ is positive definite on each $HS$.  This clearly implies that the Levi form is also positive definite, that is, $S$ is strongly pseudoconvex. A strongly pseudoconvex hypersurface will be strongly $\C$-convex if and only if $|\qf(X,X)|< \levi(X,X)$ for all non-zero $X\in HS$.

\begin{Remark}\label{R:glob}
A compact strongly $\C$-convex hypersurface $S$ also has the global property that each complex tangent hyperplane for $S$ intersects $S$ only at the point of tangency (see for example [APS, \S 2.5]).

By pushing a complex tangent hyperplane a small distance in the outward normal direction we obtain a complex hyperplane disjoint from $S$.  Moving this hyperplane to infinity by an automorphism of $\CP^n$ we see that $S$ is M\"obius-equivalent to a hypersurface in $\C^n\subsetneq \CP^n$.
\end{Remark}

Returning to the normal form in Proposition \ref{P:nrm-frm}, the strong $\C$-convexity of $S$ at $0$ is equivalent to the positive-definiteness of $\suml_{j,k=1}^{n-1} \alpha_{j,k} z_j \bar z_k +\Re \left(\suml_{j,k=1}^{n-1} \beta_{j,k} z_j z_k\right)$; this condition is equivalent in turn to the condition that the complex tangent hyperplane $\C^{n-1}\times\{0\}$ to $S$ at $0$ has minimal order of contact with $S$.

Since we may choose the constant $c$ in Proposition \ref{P:nrm-frm} to be positive, it is also clear now that $S$ is strongly $\C$-convex if and only if $S$ may be made locally strongly convex near any $p\in S$ by an automorphism of $\CP^n$.

Setting $c=0$ and diagonalizing $\qf$ with respect to a basis which is orthogonal for $\levi$ we may reduce the normal form in  Proposition \ref{P:nrm-frm} to 
\begin{equation}\label{E:exp-diag}
v= \suml_{j=1}^{n-1} \alpha_{j} |z_j|^2 +\Re \left(\suml_{j=1}^{n-1} \beta_{j} z_j^2\right) 
+ O\left(\|(z_1,\dots,z_{n-1},u)\|^3\right)
\text{ with each } \beta_j\ge 0.
 \end{equation}
Strong $\C$-convexity at the origin is now equivalent to the condition that $\beta_j<\alpha_j$ for each $j$.

\begin{Remark}\label{R:norm1}
We could normalize further to set the $\alpha_j$ equal to $1$, but we will see in \S \ref{S:dual} below that a different normalization is preferable for our purposes. 
\end{Remark}

Returning to the computations of \S \ref{S:scinv} and setting 
 \[\rho(z)=-v+\suml_{j=1}^{n-1} \alpha_{j} |z_j|^2 +\Re \left(\suml_{j=1}^{n-1} \beta_{j} z_j^2\right) 
+ O\left(\|(z_1,\dots,z_{n-1},u)\|^3\right)\] we have
\begin{align*}
\delta_{\rho,1}&=\frac{1}{4}\prod_{j=1}^{n-1} \alpha_j\\
\delta_{\rho,2}&=\frac{1}{16}
\prod_{j=1}^{n-1} \left( \alpha_j^2-\beta_j^2 \right)
\end{align*}
and so
\begin{equation}\label{E:rhoba}
\sr_S=\prod_{j=1}^{n-1}\left(1-\frac{\beta_j^2}{\alpha_j^2}\right)
\end{equation}
at $0$.

\begin{Remark}\label{R:sph1}
When $S$ is strongly $\C$-convex then $0<\sr_S\le 1$ and $\sr_S$ may be interpreted as the 
the ratio of  the volume of a sublevel set of $\levi$ to the volume of the corresponding sublevel set of $\sff$.  Moreover, the  quantity $\sr_S$ provides a scalar measure of the degree of strong $\C$-convexity; in the  case $\sr_S\equiv 1$ of ``maximal $\C$-convexity"  the hypersurface $S$ is locally M\"obius-equivalent to the unit sphere  ([Jen], [DeTr], [Bol2]).  
\end{Remark}

\subsection{Dimension two} In the next two subsections we explore additional geometric structure for strongly $\C$-convex hypersurfaces in $\CP^2$.  This material is not used in the proof of the theorems stated in \S \ref{S:bt}.

Let $p$ be a point in a strongly pseudoconvex hypersurface $S\subset\CP^2$.  By Proposition \ref{P:nrm-frm} we may assume after applying a M\"obius transformation that $p=0$ and that $S$ is given near $p$ by
\begin{equation} \label{E:exp-diag2}
v=  \alpha |z_1|^2 +\Re \left( \beta z_1^2\right) 
+ O\left(\|(z_1,u)\|^3\right)
\end{equation}
with $\alpha>0$ and $\beta\ge0$.  $S$ will be strongly $\C$-convex at $p$ if and only if $\beta/\alpha<1$.

The quantity $\beta/\alpha=\sqrt{1-\sr_S}$ defines a scalar invariant with the following interpretations.
\begin{itemize}
\item $\qf/\levi$ defines an $\R$-valued function on each $H_zS\setminus\{0\}$ which is $\R_+$-homogeneous of degree $0$ and has $\beta/\alpha$ as its maximum absolute value. 
\item Also,  the level sets of the tensor $\sff=\levi + \qf$ are ellipses in each $H_pS$ with major-to-minor axis ratio equal to $\sqrt{\frac{\alpha+\beta}{\alpha-\beta}}=\sqrt{\frac{1+\beta/\alpha}{1-\beta/\alpha}}$. 
\end{itemize}

In the two-dimensional setting we have a somewhat simpler alternate version of Proposition \ref{P:srdef} above.

\begin{Prop}\label{P:bquo} Let $S\subset\CP^2$ be a strongly pseudoconvex hypersurface, and let $R\in \Gamma(\CP^2;2,2)$ be a defining function for $S$ as in the proof of Proposition \ref{P:srdef}.
Then the quotient
\begin{equation}\label{E:bquo}
- 4
\det
\begin{pmatrix}
R_{j,k}
\end{pmatrix}
\big/
\det
\begin{pmatrix}
R_{j,\bar k}
\end{pmatrix}
\end{equation}
defines a section $\belt_S\in\Gamma(S;-3,3)$.  (Here the subscripts denote differentiation and range over $0,1,2$.)

The section $\belt_S$ is independent of the choice of defining function $R$. 

For $M\in SL(3,\C)$ we have the transformation law
\begin{equation} \label{E:btran}
\belt_S = M^* \belt_{M(S)}
\end{equation}

The corresponding scalar invariant $\left|\belt_S\right|\in\Gamma(S;0,0)$ coincides with the invariant 
$\beta/\alpha=\sqrt{1-\sr_S}$ discussed above.
\end{Prop}

\begin{proof} 
The $R_{j,k}$ are in $\Gamma(S;0,2)$, so the numerator is in $\Gamma(S;0,6)$.  The $R_{j,\bar k}$ are are in $\Gamma(S;1,1)$, so the denominator  is in $\Gamma(S;3,3)$. Thus the quotient  is in $\Gamma(S;-3,3)$.

Applying an automorphism induced by $M\in SL(3,\C)$ we find that the matrix in the numerator of \eqref{E:bquo} is multiplied on the right by $M$
and on the left by $\tensor*[^t]{M}{}$; similarly, the matrix in the denominator is multiplied on the right by $M$ and on the left by $\tensor*[^t]{\bar M}{}$. The transformation law \eqref{E:btran} follows immediately.  (Note that the lift $M^*$ is uniquely determined in this situation.)

To verify that $\belt_S$ is independent of the choice of defining function, we may procede as in the proof of Proposition 
\ref{P:srdef}; alternatively, we may use the transformation law \eqref{E:btran} along with \eqref{E:exp-diag2} above to reduce to checking at $p=(\zeta_0:0:0)$ with $R$ given by
\begin{equation*}
\tau\cdot\left(
\frac{i}{2}\zeta_0\bar{\zeta_0^2}\zeta_2
-\frac{i}{2}\zeta_0^2\bar{\zeta_0\zeta_2}
+\alpha \zeta_0\bar{\zeta_0}\zeta_1 \bar{\zeta_1}
+\frac{\beta}{2} \bar{\zeta_0^2}\zeta_1^2
+\frac{\beta}{2} \zeta_0^2\bar{\zeta_1^2} 
+\dots
\right)
\end{equation*}
where $\tau$ is  a smooth  positive section of $\oh(0,0)$ near $p$.
Using $\tau_0=\tau_{\bar 0}=0$ at $p$ we find that the expression in \eqref{E:bquo} is 
\begin{equation*}
\frac{\ \tau\beta \bar{\zeta_0^6}}
{ \tau\alpha \zeta_0^3\bar{\zeta_0^3}}
=\frac{\beta}{\alpha}\frac{\bar{\zeta_0^3}}{\zeta_0^3}
\end{equation*}
at $p$, showing  that  the quotient is  independent of the choice of defining function and also that  $\left|\belt_S\right|=\beta/\alpha$.

\end{proof}

On $\C^2$ we may use the convention of Remark \ref{R:oeuc} to write $\belt_S$ in the form
\begin{equation*}
b_S(z_1,z_2) 
\tfrac{dz_1\w dz_2}{\overset{}{d\bar{z_1}\w d\bar{z_2}}}
\end{equation*}
which describes a scalar-valued function on $\C$-linearly independent pairs $X,Y\in T_z \C^2$.  Here $b_S(z_1,z_2)=\belt_S(1,z_2,z_2)$ and  $\belt_S(\zeta_0,\zeta_1,\zeta_2)=\zeta_0^{-3}\bar{\zeta_0^3} \,b_S(\zeta_1/\zeta_0,\zeta_2/\zeta_0)$.  
We have
\begin{align*}
b_S(z_1,z_2) \tfrac{dz_1\w dz_2}{\overset{}{d\bar{z_1}\w d\bar{z_2}}}(\lam X, Y)
&=
\lam\bar{\lam\inv} \,b_S(z_1,z_2) \tfrac{dz_1\w dz_2}{\overset{}{d\bar{z_1}\w d\bar{z_2}}}(X,Y)\\
&=b_S(z_1,z_2) \tfrac{dz_1\w dz_2}{\overset{}{d\bar{z_1}\w d\bar{z_2}}}(X,\lam Y)
\end{align*}
for $\lam\in\C\setminus\{0\}$.

Writing $\rho(\zeta_0,\zeta_1,\zeta_2)=\zeta_0^2\bar{\zeta_0^2}\, r(\zeta_1/\zeta_0,\zeta_2/\zeta_0)$ we find that
\begin{equation}\label{E:aff-b}
b_S(z_1,z_2)= - \,
\det
\begin{pmatrix}
 0  & r_1  & r_2  \\
r_1  &  r_{1   1} & r_{2   1}  \\
r_{ 2 }  & r_{1   2}  &   r_{2   2} 
\end{pmatrix}
\Bigg/
\det \begin{pmatrix}
 0  & r_1  & r_2  \\
r_{\bar 1 }  &  r_{1 \bar 1} & r_{2 \bar 1}  \\
r_{\bar 2 }  & r_{1 \bar 2}  &   r_{2 \bar 2} 
\end{pmatrix}.
\end{equation}

\begin{Remark}
This type of differential is reminiscent (up to conjugation or inversion) of the Beltrami differentials $\deebar f/\dee f$ prominent in the study of quasiconformal mappings in one complex variable.
(Compare [KoRe].)
\end{Remark}

In the case where $S$ is strongly $\C$-convex (i.e., when $|\belt_S|<1$) we have seen above that $|\belt_S|=\beta/\alpha$ may be identified with the eccentricity data for a families of ellipses in each $H_z S$. The ``argument" of $\belt_S$ is similarly determined by the following condition: 

\begin{quotation}
when $X\in H_z S\setminus\{0\}$, $Y\in T_zS\setminus H_zS$ then $b_S(z_1,z_2) \tfrac{dz_1\w dz_2}{\overset{}{d\bar{z_1}\w d\bar{z_2}}}(X,Y)$ will be positive precisely when $X$ points in the direction of the minor axis of the family of ellipses in $H_zS$. \end{quotation}
(Of course, this condition becomes undefined when the ellipses are circles, i.e., when $\belt_S$ vanishes.)

\subsection{Examples} \label{SS:ex} 
\begin{enumerate}
\item The unit spheres 
\begin{equation*}
\Sigma^{(1)}_p\eqdef \{|z_1|^p + |z_2|^p=1\}
\end{equation*}
of the two-dimensional $L^p$ spaces are strongly $\C$-convex for $p>1$ and $z_1z_2\ne 0$, with
\begin{equation*}
\belt_S = \frac{2-p}{p}\frac{dz_1\w dz_2}{z_1 z_2}
\frac{\bar{z_1 z_2}}{\overset{}{d\bar{z_1}\w d\bar{z_2}}}.
\end{equation*}

\item The hypersurfaces
\begin{equation*}
\Sigma^{(2)}_\gamma\eqdef\{\Im z_2 = |z_1|^\gamma\}
\end{equation*}
are strongly $\C$-convex for $\gamma>1$ and $z_1\ne 0$, with
\begin{equation*}
\belt_S = \frac{\gamma-2}{\gamma}\frac{dz_1\w dz_2}{z_1}
\frac{\bar{z_1 }}{\overset{}{d\bar{z_1}\w d\bar{z_2}}}.
\end{equation*}

\item The hypersurfaces 
\begin{equation*}
\Sigma^{(3)}_{\alpha,\beta}\eqdef\{\Im z_2 = \alpha|z_1|^2+ \Re \beta z_1^2\}
\end{equation*}
are strongly $\C$-convex for $|\beta|<\alpha$, with
\begin{equation*}
\belt_S = \frac{\beta}{\overset{}{\alpha}}
\frac{dz_1\w dz_2}{\overset{}{d\bar{z_1}\w d\bar{z_2}}}.
\end{equation*}

\item For tube hypersurfaces $S\subset \C^2$ invariant under all real translations we have
\begin{equation*}
\belt_S = -\frac{dz_1\w dz_2}{\overset{}{d\bar{z_1}\w d\bar{z_2}}};
\end{equation*}
such hypersurfaces are never strongly $\C$-convex.

\end{enumerate}

\begin{Remark}\label{R:sph}
If $\belt_S\equiv0$ then $\sr_S\equiv1$; as noted in Remark \ref{R:sph1} this implies that $S$ is locally M\"obius-equivalent to the unit sphere.

If $\belt_S= K \frac{dz_1\w dz_2}{\overset{}{d\bar{z_1}\w d\bar{z_2}}}$, $K$ constant, $|K|\ne 0, 1,$ then $S$ is an affine image of a hypersurface of the form $\Sigma^{(3)}_{\alpha,\beta}$ above ([Bol3], [Ham]). \end{Remark}

The examples listed above have the property that
$\belt_S$ extends to a constant times a meromorphic $(2,0)$-form divided by its conjugate (reminiscent of Teichm\"uller differentials in one complex variable).  This does not hold in general.

The following result gives an indication of the restrictions that $\belt_S$ must satisfy.

\begin{Thm}\label{T:rigid}
A section
$
\lam(z_1)\,\frac{dz_1\w dz_2}{\overset{}{d\bar{z_1}\w d\bar{z_2}}}$ of $\oh(-3,3)$  
will arise as $\belt_S$
for a  strongly $\C$-convex  hypersurface
\begin{equation}\label{E:rig}
\Im z_2 = f(z_1)
\end{equation}
invariant under real $z_2$-translations if and only if
the coefficient $\lam(z_1)$ satisfies
\begin{subequations}\label{E:rig-belt-cond}
\begin{align}\label{E:rig-belt-cond-a}
\Im\left(
\lam_{\bar{z}_1\bar{z}_1}
- \bar\lam \lam_{z_1\bar{z}_1}
+\frac
{\bar\lam \lam_{\bar{z}_1}^2+\lam\lam_{\bar{z}_1}\bar{\lam}_{\bar{z}_1} - \bar{\lam^2}\lam_{z_1}\lam_{\bar{z}_1}  }
{1-\lam\bar\lam}
\right)
&=0,\\
|\lam|&<1. \label{E:rig-belt-cond-b}
\end{align}
\end{subequations}

More precisely, if $U$ is open in $\C$ and $S\subset U\times\C$ is a strongly $\C$-convex hypersurface given by \eqref{E:rig}  then the coefficient $\lam(z_1)$ of $\belt_S$ must satisfy \eqref{E:rig-belt-cond} on $U$. Conversely, when $U$ is simply-connected then any solution of \eqref{E:rig-belt-cond} gives rise to a corresponding $S\subset U\times\C$.

\end{Thm}

Note that \eqref{E:rig-belt-cond-a} may be viewed as an underdetermined hyperbolic system in the two $\R$-valued unknowns $\Re\lam$, $\Im\lam$.

\begin{Remark}  Hypersurfaces of the form \eqref{E:rig} are often known as {\em rigid} [BRT].\end{Remark}

\begin{proof}[Proof of Theorem \ref{T:rigid}]  The inequality \eqref{E:rig-belt-cond-b} is already accounted for in the definition of a strongly $\C$-convex hypersurface.

Substituting $r=  f(z_1)-\Im z_2$ into \eqref{E:aff-b} we find that the question of solving $\belt_S=\lam(z_1)\,\tfrac{dz_1\w dz_2}{\overset{}{d\bar{z_1}\w d\bar{z_2}}}$ with $S$ as in \eqref{E:rig} reduces to that of the solvability of 
\begin{equation*}\label{E:pre-2be-rig}
f_{z_1 z_1}=\lam(z_1) f_{z_1\bar{z}_1} 
\end{equation*}
for $\R$-valued $f$.
For typographical simplicity we drop  subscripts to rewrite this as
\begin{subequations}\label{E:2be-rig}
\begin{align} f_{z z}&=\lam(z) f_{z\bar{z}}
\label{E:2be-rig-a}\\
\bar f&= f.\label{E:2be-rig-b}
\end{align}
\end{subequations}

Differentiating \eqref{E:2be-rig-a} with respect to $\bar z$ we get
\begin{equation*}
f_{zz\bar z} = \lam f_{z\bar{zz}} + \lam_{\bar z} f_{z\bar z}.
\end{equation*}

Conjugating and applying \eqref{E:2be-rig-b} we have
\begin{equation*}
f_{z\bar{zz}} = \bar\lam f_{zz\bar{z}} + \bar\lam_{z} f_{z\bar z}.
\end{equation*}

Since we are assuming in particular that $S$ is strongly pseudoconvex we have $f_{z\bar z}>0$. Setting $h=\log f_{z\bar z}$ we have
\begin{subequations}\label{E:2be-rig-prev}
\begin{align}
h_z &= \lam h_{\bar z}+\lam_{\bar z}\\
h_{\bar z} &= \bar \lam h_{ z}+\bar\lam_{ z}\\
\bar h &= h.
\end{align}
\end{subequations}

Using linear algebra to isolate $h_z$ and $h_{\bar z}$ this may be rewritten as
\begin{subequations}\label{E:2be-rig-rev}
\begin{align}
h_z &=\frac{\lam_{\bar z} + \lam \bar\lam_z}{1-\lam\bar \lam}\label{E:2be-rig-rev-a}\\
h_{\bar z} &= \frac{\bar\lam \lam_{\bar z} +\bar\lam_z}{1-\lam\bar \lam}\label{E:2be-rig-rev-b}\\
\bar h &= h.\label{E:2be-rig-rev-c}
\end{align}
\end{subequations}

Differentiating \eqref{E:2be-rig-rev-a} with respect to $\bar z$ and matching this with the result of differentiating \eqref{E:2be-rig-rev-b} with respect to $z$ we obtain \eqref{E:rig-belt-cond-a}.

For the converse, note that from the previous paragraph we see that \eqref{E:rig-belt-cond-a} is precisely the condition guaranteeing that the form
\begin{equation}\label{E:form}
\frac{\lam_{\bar z} + \lam \bar\lam_z}{1-\lam\bar \lam}\,dz
+
\frac{\bar\lam \lam_{\bar z} +\bar\lam_z}{1-\lam\bar \lam}\,d\bar z
\end{equation}
is $d$-closed.  Since \eqref{E:form} is self-conjugate we see that \eqref{E:rig-belt-cond-a} is precisely the condition required to solve \eqref{E:2be-rig-rev} on  simply-connected $U$.  

It remains to show that solutions of \eqref{E:2be-rig-rev} (equivalently, of \eqref{E:2be-rig-prev}) give rise to solutions of \eqref{E:2be-rig}.
We begin by solving $g_{z\bar z}=e^h,$ $\bar g=g$ on $U$.  We have 
\begin{equation*}
\left(
g_{z z}-\lam g_{z\bar z}
\right)_{\bar z}
= e^h
\left(
h_z - \lam_{\bar z} - \lam h_{\bar z}\right)=0,
\end{equation*}
so we may write
\begin{equation*}
g_{z z}=\lam g_{z\bar z}- H_{zz}
\end{equation*}
with $H$ holomorphic.  Writing
\begin{equation*}
f=g+H+\bar H
\end{equation*}
we have $\bar f=f$ and $f_{zz}=g_{zz}+H_{zz} = \lam g_{z\bar z} = \lam f_{z\bar z}
$
as required.
\end{proof}

\begin{Q}
What conditions must $\belt_S$ satisfy for general (strongly $\C$-convex) $S$?
\end{Q}

\section{Dual hypersurfaces}\label{S:dual}

Let $\CP^{n*}$ denote the projective space dual to $\CP^n$.  Each point $\zeta^*=(\zeta_0^*:\ldots:\zeta_n^*)$ in $\CP^{n*}$ determines a hyperplane
\begin{equation*}
\hyp_{\zeta^*}\eqdef
\{\zeta\in\CP^n \st \zeta_0 \zeta_0^* + \dots +\zeta_n\zeta_n^*=0\}
\end{equation*}
in $\CP^n$; conversely, for $\zeta\in\CP^n$ the set
$\{\zeta^*\in \CP^{n*} \st \zeta\in\hyp_{\zeta^*}\}$ defines a hyperplane $\hyp_\zeta^*$ in $\CP^{n*}$. 

 We define the {\em incidence manifold} $\inc\subset\CP^n\times\CP^{n*}$ by
 \begin{align*}
\inc &=
\{ (\zeta,\zeta^*)\in \CP^n\times\CP^{n*}
\st \zeta_0 \zeta_0^* + \dots +\zeta_n\zeta_n^*=0\}\\
&= \{ (\zeta,\zeta^*)\in \CP^n\times\CP^{n*}
\st \zeta\in \hyp_{\zeta^*}\}\\
&= \{ (\zeta,\zeta^*)\in \CP^n\times\CP^{n*}
\st \zeta^*\in \hyp_{\zeta}\}.
\end{align*}
$\inc$ is a complex manifold of dimension $2n-1$.

For $M\in SL(n+1,\C)$ the automorphism
\begin{align}\label{E:dual-aut}
\mapdef{\psi_M: \CP^n\times\CP^{n*}}
{\CP^n\times\CP^{n*}}
{(\zeta,\zeta^*)}
{(M\zeta,\tensor*[^t]{M}{^{-1}}\zeta^*)}\notag
\end{align}
maps $\inc$ to $\inc$.  The action of $SL(n+1,\C)$ on $\inc$ is transitive.

A smooth real hypersurface $S\subset \CP^n$ induces a map $\dual_S:S\to\CP^{n*}$ defined by the rule that $\hyp_{\dual_S(\zeta)}$ is the unique complex hyperplane in $\CP^n$  tangent to $S$ at $\zeta$. (For $S\subset \C^n\subsetneq\CP^n$ we may write this in affine coordinates as $\hyp_{\dual_S(z)}=z+H_z(S)$.) 

Let $\gr_S\subset \CP^n\times\CP^{n*}$ denote the graph of $\dual_S$.  Then $\gr_S\subset\inc$ and we have the transformation law $\gr_{M(S)}=\psi_M(\gr_S)$.

\begin{Thm}
When $S$ is compact and strongly $\C$-convex then $\dual_S$ maps $S$ diffeomorphically onto a smooth strongly $\C$-convex hypersurface $S^*\subset\CP^{n*}$.  This map satisfies the contact condition
\begin{equation}\label{E:ctact}
\dual'_S(\zeta)\left(H_\zeta S\right)=H_{\dual_S(\zeta)}S^*
\end{equation}
but this map is never $\C$-linear and thus $\dual_S$ is never a $CR$ map.

We have $S^{**}=S$ and $\dual_{S^*}\circ \dual_S=I_S$.

Furthermore
\begin{equation}\label{E:sr-dual}
\sr_{S^*}\circ \dual_S = \sr_{S},
\end{equation}
hence in the special case $n=2$ 
\begin{equation}\label{E:belt-dual}
|\belt_{S^*}|\circ \dual_S = |\belt_{S}|.
\end{equation}
\end{Thm}

Little or none of this is new, but for convenience we provide below a presentation of the purely local parts of this result.  For the global aspects see [APS, 
\S 2.5] (as well as [MT]).

\subsection{Examples}  Returning to the examples of \S \ref{SS:ex}, it is not hard to verify that the dual of $\Sigma^{(1)}_{p}$ is M\"obius-equivalent to the standard dual 
 $\Sigma^{(1)}_{p/(p-1)}$;  that the dual of $\Sigma^{(2)}_{\gamma}$ is M\"obius-equivalent to $\Sigma^{(2)}_{\gamma/(\gamma-1)}$; and that $\Sigma^{(2)}_{\alpha,\beta}$ is self-dual (up to M\"obius equivalence).
 
 \subsection{Affine coordinates}\label{S:affco} It will be convenient to use affine coordinates $z$ for $\CP^n$ and $\eta$ for
 $\CP^{n*}$ given by 
 \begin{align*}
(z_1,\dots,z_n)&=(\zeta_1/\zeta_0,\dots,\zeta_n/\zeta_0)\\
(\zeta_0:\zeta_1:\dots:\zeta_n)&=(1:z_1:\dots:z_n)\\
(\eta_1,\dots,\eta_{n-1},\eta_n)
&=(-\zeta^*_1/\zeta^*_n,\dots,-\zeta^*_{n-1}/\zeta^*_n,\,\zeta^*_0/\zeta^*_n)\\
(\zeta^*_0:\zeta^*_1:\dots:\zeta^*_{n-1}:\zeta^*_n)&=
(\eta_n:-\eta_1:\dots:-\eta_{n-1}:1).
\end{align*}

In $(z,\eta)$-coordinates the incidence manifold $\inc$ is given by
\begin{equation*}
z_n+\eta_n=z_1\eta_1+\dots+z_{n-1}\eta_{n-1}.
\end{equation*}

This choice of affine coordinates has good compatibility with Proposition \ref{P:nrm-frm} and its variants in that
the point $(1:0:\dots:0)$ with affine coordinates $z=0$ is mapped by $\dual_S$ to the point $(0:\dots:0:1)$ with affine coordinates $\eta=0$.

\subsection{Dimension two}  To study $\dual_S$ locally near a point $\zeta\in S\subset\CP^2$ we may first apply an automorphism of $\CP^2$ to reduce to the case where $\zeta=(1:0:0)$ and $S$ is given locally by \eqref{E:exp-diag2} with respect to affine coordinates $(z_1,u+iv)$.  Writing $\eta=\dual_S(z)$ using affine coordinates as above we see that the complex line tangent to $S$ at $z$ is given by
\begin{subequations}\label{E:leg}
\begin{align}
z_2+\eta_2&=z_1\eta_1.\label{E:leg-a}\\
\intertext{
Since this line is parallel in $\C^2$ to $H_zS$ we have
}
dz_2&= \eta_1\,dz_1 \text{ on } H_zS.\label{E:leg-b}
\end{align}
\end{subequations}

For $z\in S$ near zero, \eqref{E:leg-b} determines $\eta_1$ and then \eqref{E:leg-a} determines $\eta_2$. From this and \eqref{E:exp-diag2} we deduce that the derivative of $\dual_S$ at $0$ 
takes the form
\begin{align}\label{E:dsp}
\mapdef{\dual_S'(0):T_0S=\C\times\R}
{T_0 S^*=\C\times\R\notag}
{\begin{pmatrix} z_1 \\ u \end{pmatrix}}
{\begin{pmatrix} 2i\beta z_1 + 2i\alpha\bar z_1 \\ -u \end{pmatrix}}.
\end{align}
In particular we see that $\dual_S'(0)$ maps $H_0S=\C\times\{0\}$ to  $H_0S^*=\C\times\{0\}$ but the assumption that $\alpha\ne 0$ guarantees that this  map is not $\C$-linear; it follows that $\dual_S$ is contact but not $CR$.   Differentiating \eqref{E:leg-a} along $\inc$ and subtracting \eqref{E:leg-b} using the contact condition \eqref{E:ctact} we find that
\begin{equation}\tag{\ref{E:leg}c}
d\eta_2=z_1\,d\eta_1 \text{ on } H_\eta S^*.
\end{equation}

The equations (\ref{E:leg}c) and \eqref{E:leg-a} allow us to determine $z\in S$ from $\eta\in S^*$.  The symmetry of the equations \eqref{E:leg} shows that 
$\dual_{S^*}\circ \dual_S=I_S$ and hence that $S^{**}=S$.

\begin{Remark}
$\dual'_S(0)$ maps the ellipses in $H_0S$ determined by $\sff_S$ to the ellipses in $H_0S^*$ determined by $\sff_{S^*}$, switching major and minor axes. 
\end{Remark} 

\begin{Remark}
The equations \eqref{E:leg} are reminiscent of the equations describing the Legendre transform of a strongly convex real planar curve (see for example [H\"or, p. 18]).
\end{Remark}

Using the description of $\dual_S'(0)$ above 
we find that 
\begin{align*}
\eta_2 &= z_1\eta_1- z_2\\
&= z_1 \left(2i \beta z_1 + 2i\alpha \bar z_1  
+ O\left(\|(z_1,u)\|^2\right)  \right)
-u-i\left( \alpha |z_1|^2+\Re \beta z_1^2 +  O\left(\|(z_1,u)\|^3\right)  \right).
\end{align*}
Combining this with
\begin{equation*}
z_1 = \frac{i\beta\eta_1+i\alpha \bar\eta_1}
{2(\alpha^2-\beta^2)}
+O\left(\|(\eta_1,\Re\eta_2)\|^2\right)
\end{equation*}
we obtain
\begin{align}
\Im \eta_2 &= \Im \left(z_1\eta_1- z_2\right)\notag\\
&=  \alpha |z_1|^2+\Re \beta z_1^2 +  O\left(\|(z_1,u)\|^3\right)\label{E:exp-dual2}\\
&= \frac{\alpha|\eta_1|^2+\Re\beta\eta_1^2}
{4(\alpha^2-\beta^2)}+O\left(\|(\eta_1,\Re\eta_2)\|^3\right)\notag
\end{align}
verifying \eqref{E:sr-dual} and \eqref{E:belt-dual}.

Referring back to Remark \ref{R:norm1} we see that the normalization $\alpha=1$ would not carry over to the dual hypersurface.   But we may apply coordinate dilations to arrange instead that 
\begin{equation}\label{E:4n1}
\alpha^2-\beta^2=1/4,
\end{equation}
and we see from \eqref{E:exp-dual2} that this normalization does carry over to the dual.

\subsection{Higher dimension}\label{S:dual-hd} In higher dimension we may argue as above, replacing the equations  \eqref{E:leg}
 by
\begin{subequations}\label{E:leg-hd}
\begin{align}
z_n+\eta_n&=z_1\eta_1+\dots+z_{n-1}\eta_{n-1}\label{E:leg-hd-a}\\
dz_n&= \eta_1\,dz_1+\dots+\eta_{n-1}\,dz_{n-1} \text{ on } H_zS \\
d\eta_n&=  z_1\,d\eta_1+\dots+z_{n-1}\,d\eta_{n-1} \text{ on } H_\eta S^*.
\end{align}
\end{subequations}

To study the duality near a general point of $S$ we may normalize as before to reduce to the study of \eqref{E:exp-diag} 
near $0$.  Following the lead of the previous section, we may use coordinate dilations to further assume that
\begin{equation}\label{E:4n}
\alpha_j^2-\beta_j^2=1/4.
\end{equation}
Then \eqref{E:dsp} is replaced by 
\begin{align}\label{E:dsp-hd}
\mapdef{\dual_S'(0):T_0S=\C^{n-1}\times\R}
{T_0 S^*=\C^{n-1}\times\R\notag}
{\begin{pmatrix} z_1 \\ \vdots \\ z_{n-1} \\ u \end{pmatrix}}
{\begin{pmatrix} 2i \beta_{1}z_1+2i\alpha_1\bar z_1  \\
\vdots\\
2i\beta_{n-1}z_{n-1}+2i\alpha_{n-1}\bar z_{n-1}\\ 
 -u \end{pmatrix}}
\end{align}
and \eqref{E:exp-dual2}
is replaced by
\begin{equation}\label{E:exp-dual-n}
 \Im \eta_n = \suml_{j=1}^{n-1} \alpha_j |\eta_j|^2  +\Re \left(\suml_{j=1}^{n-1} \beta_j \eta_j^2\right)
+ O\left(\|(\eta_1,\dots,\eta_{n-1},\Re \eta_n)\|^3\right).
\end{equation}
The transformation law \eqref{E:sr-dual} now follows from \eqref{E:rhoba}.

Note that \eqref{E:dsp-hd} shows that $\dual_S$ is orientation-preserving when $n$ is even and orientation-preserving when $n$ is odd (when both $S$ and $S^*$ are oriented as the boundary of the corresponding pseudoconvex domain).

We further note for future reference that the normalization 
\eqref{E:4n} allows us to rewrite \eqref{E:rhoba} as
\begin{equation}\label{E:rhoban}
\sr_S=\prod_{j=1}^{n-1}\frac{1}{4\alpha_j^2}.
\end{equation}

\begin{Remark} \label{R:anti-cr}
From \eqref{E:dsp-hd} we see that the maps \eqref{E:ctact} will be anti-$\C$-linear for all $\zeta\in S$ if and only if the $\beta_j$ always vanish, that is, if and only if $\sr_S\equiv1$.  From Remark \ref{R:sph1}, it follows that $\dual_S$ is anti-$CR$ if and only if $S$ is locally a M\"obius image of a sphere. 
\end{Remark}

\subsection{Dimension one} \label{S:dual-1}  When $n=1$ it is often preferable to replace the affine coordinates from \S\ref{S:affco}  by 
\begin{align*}
z&=\zeta_1/\zeta_0\\
(\zeta_0:\zeta_1)&=(1:z)\\
\eta
&=-\zeta^*_{0}/\zeta^*_1\\
(\zeta^*_0:\zeta^*_1)&=
(-\eta:1).
\end{align*}
In these coordinates, the incidence manifold $\inc$ is described by $z=\eta$, the map $\dual_S$ is the identity, and $S^*=S$.

\section{The bilinear pairing}\label{S:bil}
The expression
\begin{equation*}
\frac{(-2)^{1-n}}{(n-1)!} \suml_{j=0}^n \zeta^*_j \,d \zeta_j 
\w \left( \suml_{j=0}^n d\zeta_j \w d \zeta^*_j \right)^{n-1}
\end{equation*}
defines a holomorphic $\oh(n,0)\times\oh(n,0)$-valued $(2n-1,0)$-form on $\CP^n\times\CP^{n*}$.  Let $\nu$ denote the restriction of this form to the incidence manifold $\inc$.  Since $\dim_\C\inc=2n-1$, $\nu$ is closed along $\inc$.

Since 
\begin{equation*}
\suml_{j=0}^n \zeta^*_j \,d \zeta_j 
+ \suml_{j=0}^n \zeta_j \,d \zeta^*_j 
=d\left( \suml_{j=0}^n \zeta_j  \zeta^*_j \right) = 0 
\end{equation*}
on $\inc$, $\nu$ is symmetric with respect to interchange of $\zeta$ with $\zeta^*$ when $n$ is even and  antisymmetric when $n$ is odd.  

For $M\in SL(n+1,\C)$ it is easy to check that $\psi_M^* \nu=\nu$, where $\psi_M$ was defined in \eqref{E:dual-aut}.  From the $SL(n+1,\C)$-homogeneity of $\inc$ we see that $\nu$ is determined up to a multiplicative constant by this invariance property.  

Let $S\subset\CP^n$ be compact and strongly $\C$-convex, and let $F\in\Gamma(S;-n,0)$ and  $G\in\Gamma(S^*;-n,0)$.  Then $FG\,\nu$ defines a standard $(2n-1)$-form on $\inc$, so we may define 
\begin{equation}\label{E:bil-def}
\lma F,G \rma=\lma F,G \rma_{S,S^*} \eqdef \intl_{\Gamma_S} FG\,\nu
\end{equation}
as the pairing discussed in Theorem 
\ref{T:bt}.
Morever, $FG\,\nu$ will be holomorphic and hence closed in any portion of $\inc$ where $F$ and $G$ are holomorphic; item \itemref{N:bt}{I:bt-iso} from Theorem \ref{T:bt} follows immediately.
Item \itemref{N:bt}{I:bt-mob} of Theorem \ref{T:bt} follows from the invariance  of $\nu$ described in the previous paragraph.

For $\zeta=\{(1:0:\dots:0)\}, \zeta^*=\{(0:\dots:0:1)\}$ we may use the conventions of \S\ref{S:affco} and Remark \ref{R:oeuc} to write 
\begin{equation*}
\nu = 2^{1-n} \frac
{dz_1 \w d\eta_1\w\dots\w dz_{n-1} \w d\eta_{n-1} \w dz_n}
{\sqrt[n+1]{dz_1\w\dots\w dz_n}\sqrt[n+1]{d\eta_1\w\dots\w d\eta_n}} .
\end{equation*}
Writing $F=f(z) \sqrt[n+1]{dz_1\w\dots\w dz_n}$, $G=g(\eta)\sqrt[n+1]{d\eta_1\w\dots\w d\eta_n}$ we see that the integrand $FG\,\nu$ may be written as \[2^{1-n}f(0) g(0) \,
dz_1 \w d\eta_1\w\dots\w dz_{n-1} \w d\eta_{n-1} \w dz_n\] at the point in question.
Assuming now the normalizations \eqref{E:exp-diag} and \eqref{E:4n} we see from \eqref{E:dsp-hd} that 
\begin{equation*}
d\eta_j = 2i\beta_j z_j + 2i \alpha_j\bar z_j
\end{equation*}
on $T_{(0,0)}\Gamma_S$; recalling \eqref{E:rhoban} this allows us to further rewrite the integrand  as
\begin{align}\label{E:rewrite-bil}
FG\,\nu &=(-2i)^{1-n}f(0) g(0) \phi_S(0)^{-1/2}\,
dz_1 \w d\bar z_1\w\dots\w dz_{n-1} \w d\bar z_{n-1} \w dz_n \notag\\
&= f(0) g(0) \phi_S(0)^{-1/2}\,dS 
\end{align}
at $z=0$, 
where
\begin{equation*}
dS=dx_1\w dy_1\w\dots\w dx_{n-1}\w dy_{n-1}\w dx_n
\end{equation*}
is the element of euclidean surface area for $S$ at $0$.

Writing everything in terms of $\eta$ variables instead we may also write 
\begin{equation}\label{E:rewrite-bil-dual}
FG\,\nu = (-1)^n f(0) g(0) \phi_S(0)^{-1/2}\,dS^*
\end{equation}
at $\eta=0$.

\section{An invariant Hardy space}\label{S:hard}

We begin by constructing an invariantly-defined positive $\oh(n,n)$-valued $(2n-1)$-form $\mu=\mu_S$ on $S$.  This will allow us to set
\begin{equation}\label{E:norm-def}
\|F\|_S^2 = \intl_S F\bar F\,\mu
\end{equation}
for $F\in\Gamma(S: j,k)$ with $j+k=-n$.

There is such a form already in the literature, namely {\em Fefferman's measure} $\muf$ given in affine (or rather, euclidean) terms by
\begin{equation}\label{E:Fef-def}
  \frac{2\left|
\det\begin{pmatrix}
0 & \rho_{j}\\
\rho_{\bar k} & \rho_{j \bar k}
\end{pmatrix}
\right|^{\frac{1}{n+1}}
\,dS}{\|d\rho\|
  \left(dz_1\w\dots\w dz_n\right)^{n/(n+1)}\left(d\bar z_1\w\dots\w d\bar z_n\right)^{n/(n+1)}},
\end{equation}
 where $\rho$ is a defining function for $S$ and $dS$ is the euclidean surface area form on $S$
([Fef, p. 259]; see also [Bar]).

\begin{Remark} Fefferman's definition included an unspecified dimensional constant. In [Bar] a preferred choice of dimensional constant is set forth,  but a different choice has been made here for simplicity in the current context. The choice of this dimensional constant has no consequences outside of the current section.\end{Remark}

The form $\muf$ is invariant not only under M\"obius tranformations but under general $CR$ diffeomorphism (in any context where the roots of the canonical bundle can be handled properly).  Unfortunately, as we shall see in a moment, $\muf$ does not connect optimally with the pairing constructed in \S \ref{S:bil}.  But the form $\mu$ we seek must differ from $\muf$ by multiplication by some scalar M\"obius-invariant quantity.

Picking a point $p\in S$ and normalizing as before using 
\eqref{E:exp-diag} and \eqref{E:4n} we have
\begin{equation*}
\muf = \frac{\sr_S^{-1/(2n+2)}\, dS}
{\left(dz_1\w\dots\w dz_n\right)^{n/(n+1)}\left(d\bar z_1\w\dots\w d\bar z_n\right)^{n/(n+1)}}
\end{equation*}
at $0$.  But referring back to \eqref{E:rewrite-bil} it is clear that we would rather have
\begin{equation}\label{E:mu-predef}
\mu = \frac{\sr_S^{-1/2}\, dS}{\left(dz_1\w\dots\w dz_n\right)^{n/(n+1)}\left(d\bar z_1\w\dots\w d\bar z_n\right)^{n/(n+1)}}
 = \sr_S^{-n/(2n+2)} \muf.
\end{equation}
So we define $\mu$ to be $\sr_S^{-n/(2n+2)} \muf$.

The M\"obius-invariance of $\mu$ follows from that of 
$\sr_S$ and  $\muf$.  For a direct proof of this fact, note that it suffices to check that $\mu$ is invariant at $0$ under maps of the type described in Remark \ref{R:iso}. But under such a map the numerator in the middle of \eqref{E:mu-predef} picks up a factor of $|\det M'|^2 \kappa=\kappa^n$ whereas the denominator picks up a factor of $\left| \det M'\cdot \kappa\right|^{2n/(n+1)}=\kappa^n$; the cancellation proves the invariance.

Item \itemref{N:an}{I:aff-norm-def} of Theorem \ref{T:aff-norm} follows now from  \eqref{E:Fef-def} and 
Proposition \ref{P:srdef} above.

Using \eqref{E:mu-predef}, \eqref{E:rewrite-bil} and \eqref{E:rewrite-bil-dual} we obtain the pointwise relation
\begin{equation*}
\left| FG\,\nu\right| = |F| \sqrt{\mu_S} \cdot |G| \sqrt{\mu_{S^*}};
\end{equation*}
upon integration we obtain   Cauchy-Schwarz inequality 
\begin{equation}\label{E:CS}
\left\|\lma F, G \rma\right\|
\le \| F \|_S\| G \|_{S^*}.
\end{equation}
Item \itemref{N:bt}{I:bt-xd} of Theorem \ref{T:bt} follows now upon observing that equality holds in \eqref{E:CS} when $\arg FG\,\nu$ is constant and the normalized coefficients of $F$ and $G$ agree in absolute value.

\begin{Remark} \label{R:uniq}
Any invariant norm satisfying \itemref{N:bt}{I:bt-xd} would have to formed by multiplying $\mu$ above by a positive scalar invariant $\chi_S$ satisfying $\chi_S\cdot \left(\chi_{S^*}\circ \dual_S\right)=1.$  But if the norm is to be based on second-order information then from \eqref{E:exp-dual-n} we see that we must also have $\chi_S=\chi_{S^*}\circ \dual_S$.  Hence $\chi_S=1$ showing that our pairing is uniquely determined by the given conditions.
\end{Remark}

\begin{proof}[Proof of Theorem \ref{T:aff-norm}, item \itemref{N:an}{I:no-norm}]
Any such norm would come from an invariantly-defined positive  $\oh(j+k,j+k)$-valued $(2n-1)$-form $\tilde\mu=\tilde\mu_S$ on $S$.  Consider the special case of the hyperquadric
\begin{equation*}
S=\left\{v= \frac12 \suml_{j=1}^{n-1}  |z_j|^2\right\}.
\end{equation*}
We may write
\begin{equation*}
\tilde\mu=
\frac{u\, dS}
{\left(dz_1\w\dots\w dz_n\right)^{(j+k)/(n+1)}\left(d\bar z_1\w\dots\w d\bar z_n\right)^{(j+k)/(n+1)}}
\end{equation*}
where $u$ is a positive scalar function on $S$.  Invariance of $\tilde\mu$ under the non-isotropic dilations
\begin{equation*}
(z_1, \dots, z_{n-1},z_n)\mapsto
(\lam z_1, \dots, \lam z_{n-1},\lam^2 z_n)
\end{equation*}
forces $\frac{\lam^{2n} u(0)}{\lam^{2(j+k)}}=u(0)$ and hence $u(0)=0$, contradicting the positivity of $\tilde\mu$.
\end{proof}

\begin{Remark}
The norm constructed above appeared earlier in [BaLa] in the special setting  of two-dimensional Reinhardt domains, where it is obtained by integration with respect to the measure $\mu_0$ set forth in \S 8 of that work.
\end{Remark}

\section{The invariant projection operator} \label{S:ler}
We let
\begin{equation*}
\Phi(\zeta,\zeta^*) = \left( \suml_{j=0}^n \zeta_j \zeta^*_j\right)^{-n};
\end{equation*}
$\Phi$ may be viewed as a meromorphic section of $\oh(-n,0)\times \oh(-n,0)$  on $\CP^n\times\CP^{n*}$ with pole along $\inc.$
We may use $\Phi$ to   define \begin{equation*}
\Phi_\zeta:\zeta^*\mapsto \Phi(\zeta,\zeta^*).
\end{equation*}

\begin{Thm} \label{T:ler-lerc} If $S$ be a compact strongly $\C$-convex real hypersurface in $\CP^n$ then
the formula
\begin{equation}\label{E:lerayd}
(\ler_S F)(\zeta)
=\frac12 F(\zeta)
+\frac{(n-1)!}{2}\left(\frac{i}{\pi}\right)^n\cdot\pv\lma F, \Phi_\zeta\rma
\end{equation}
defines a bounded projection operator
\begin{equation}\label{E:lerayspec}
\ler_S:  L^2(S;-n,0)
\to \hs (S;-n,0).
\end{equation}
(Here $\pv$ denotes the principal value of the singular integral.)

For any lift  $M^*$ (as in \S \ref{S:ojk}) of an automorphism of $\CP^n$ mapping $S_1$ to $S_2$ we have 
\begin{equation}\label{E:lerayinv}
M^*(\ler_{S_2}F)=\ler_{S_1}(M^* F).
\end{equation}
\end{Thm}

\begin{Remark}
For $n=1$ using the affinizations from \S  \ref{S:dual-hd}   it is easy to check that $\ler_S=\cau_+$ as defined in \ref{E:CauchyTransform}.
\end{Remark}

\begin{proof}[Proof of Theorem \ref{T:ler-lerc}] The invariance is clear from the construction.

Using the invariance together with Remark \ref{R:glob} we may assume that $S\subset\C^n$. Note also that Remark \ref{R:glob} also guarantees that the principal value operator above is singular only along the diagonal.

We claim that after converting to affine coordinates as in \S \ref {S:dual-hd} and using the standard trivializations of $\oh(-n,0)$ over $\C^n$, the operator $\ler_S$
coincides now with the classic {\em Leray transform}
$\lerc_S$ ([Ler], see also [A\u{i}z1]) defined by
\begin{equation}\label{E:ler-cl}
\left(\lerc_S f\right)(z)=
\frac12 f(z) + (2\pi i)^{-n} \pv  \int\limits_{w\in S}
f(w)
\frac{\dee\rho(w)\w\left(\deebar\dee\rho(w)\right)^{n-1}}
{\left(\dee\rho(w)[w-z]\right)^n}
\end{equation}
where $\dee\rho(w)[w-z]=\suml_{j=1}^n \frac{\dee\rho}{\dee w_j}(w)(w_j-z_j)$ and $\pv$ again denotes the principal value of the singular integral;
then the remaining claims follow from well-known facts about $\lerc_S$ [KeSt1] (see also [Han]).

It suffices to show that both integral terms match for $z,w\in S$.  Since $\lerc_S$ is known to have the same invariance property noted above for $\ler_S$ [Bol2] it will suffice to check this under the assumption that $w=0$ and $S$ is given near $0$ by \eqref{E:exp-diag}  with further normalization \eqref{E:4n}.  Then routine computation reveals that in either formulation the contribution to the integral term at $w$ is
\begin{equation*}
\frac{(n-1)!f(0)}{4} \left( \frac{2i}{\pi  z_n} \right)^n 
\left(\prod_{j=1}^{n-1} \alpha_{j}\right) \,dS
\end{equation*}
where $dS$ is euclidean surface area.
\end{proof}

\begin{Thm}\label{T:sal}
For $S$ as above we have
\begin{equation}\label{E:sal}
\lma \ler_S F, G \rma = \lma F, \ler_{S^*} G \rma
=\lma \ler_S F, \ler_{S^*} G \rma.
\end{equation}
\end{Thm}

See the proof of Theorem 26 in [Lin] for closely related facts.

\begin{proof}
The first equality follows from \eqref{E:lerayd} and Fubini's theorem.  (To accommodate the use of singular integrals here, we do this
for a sequence of modified versions of $\Phi$ truncated near the singularity, then pass to the limit.) Then we also have
$\lma \ler_S F, G \rma =
\lma \ler_S^2 F, G \rma =
\lma \ler_S F, \ler_{S^*} G \rma.$
\end{proof}

The conditions \eqref{E:lerayspec} and \eqref{E:sal} characterize $\ler_S$ (since a projection operator is determined by its image and kernel).

The following three corollaries complete the proof of Theorem \ref{T:bt}.

\begin{Cor} \label{C:kee}
In the above setting we have
\begin{align}\label{E:lerayn}
\left\|\ler_S\right\|{}
&=
\sup _{\substack{
 F\in L^2(S;-n,0), \,\left\|  F \right\|{}\le 1\\
G\in L^2(S^*;-n,0), \,\left\| G \right\|{}\le 1
}}
\left|\lma \ler_S  F, G \rma\right|\notag\\
&=\sup _{\substack{
 F\in L^2(S;-n,0), \,\left\|  F \right\|{}\le 1\\
G\in L^2(S^*;-n,0), \,\left\| G \right\|{}\le 1
}}
\left|\lma  F, \ler_{S^*} G \rma\right|\notag\\
&=
\sup _{\substack{
 F\in L^2(S;-n,0), \,\left\|  F \right\|{}\le 1\\
G\in L^2(S^*;-n,0), \,\left\| G \right\|{}\le 1
}}
\left|\lma \ler_S  F, \ler_{S^*} G \rma\right|\\
&=\left\|\ler_{S^*}\right\|{}. \notag
\end{align}
\end{Cor}

\begin{proof}
This follows from item  \itemref{N:bt}{I:bt-xd}  of Theorem \ref{T:bt} along with Theorem \ref{T:sal}.
\end{proof}

\begin{Cor}\label{C:test}
For $F\in \hs(S)$ we have
\begin{align*}
\frac{\left\| F\right\|{}}{\left\|\ler_S\right\|{}}
\le \sup_{h\in\hs(S^*), \left\| h \right\|{} \le 1}
\left|\lma F,h \rma\right|
\le \left\| F\right\|{}.
\end{align*}
\end{Cor}

\begin{proof}
The right-hand inequality follows from \eqref{E:CS}.

For the left-hand inequality we cite item  \itemref{N:bt}{I:bt-xd}  of Theorem \ref{T:bt}, Theorem \ref{T:ler-lerc} and Theorem \ref{T:sal} to obtain
\begin{align*}
\left\| F\right\|{}
& = \sup_{G\in L^2(S^*;-2,0),\left\| G\right\|{} \le 1} \left|\lma F, G \rma\right\|\\
&=  \sup_{G\in L^2(S^*;-2,0),\left\| G\right\|{} \le 1} \left|\lma \ler_S F, G \rma\right|\\
& = \sup_{G\in L^2(S^*;-2,0),\left\| G\right\|{} \le 1} \left|\lma F, \ler_{S^*} G \rma\right|\\
& \le \sup_{H\in\hs(S^*),\left\| H\right\|{} \le \left\|\ler_S\right\|{} } \left|\lma F, H \rma\right|\\
&= \left\|\ler_S\right\|{}  \sup_{H\in\hs(S^*),\left\| H\right\|{} \le 1} \left|\lma F, H \rma\right|.
\end{align*}
\end{proof}

\begin{Remark}
Corollary \ref{C:test} may be viewed as a Hardy space version of the duality theorem of Martineau [Mar] and A\u{i}zenberg [A\u{i}z2] -- see also [APS, Chapter 3] and especially [Lin, Thm. 26].
\end{Remark}

\begin{Cor}\label{C:pc}
In the above setting we have
\begin{equation}\label{E:canc}
\inf_{F\in \hs(S), \left\| F\right\|{}=1}
\sup_{H\in \hs(S^*), \left\| H\right\|{}\le 1}
\left|\lma F,H \rma\right|
=\frac{1}{\left\|\ler_S\right\|{}}.
\end{equation}
\end{Cor}

\begin{proof}
Corollary \ref{C:test} shows that the left-hand side of \eqref{E:canc} is $\ge \frac{1}{\left\|\ler_S\right\|{} }$.  

For the other half we note that for small $\eps>0$ we may pick $\tilde F\in L^2(S;-n,0)$ and $F=\ler_S \tilde F \in \hs(S)$ with $\left\| F\right\|{} =1$, $\left\| \tilde F\right\|{} \le \frac{1}{\left\|\ler_S\right\|{}-\eps}.$ Then 
\begin{align*}
\sup_{H\in \hs(S^*), \left\| H\right\|{}\le 1}
\left|\lma F,H \rma\right|
&=
\sup_{H\in \hs(S^*), \left\| H\right\|{}\le 1}
\left|\lma \ler_S \tilde F,H \rma\right|\\
&=
\sup_{H\in \hs(S^*), \left\| H\right\|{}\le 1}
\left|\lma  \tilde F,H \rma\right|\\
&\le \left\| \tilde F \right\|{}\\
&\le \frac{1}{\left\|\ler_S\right\|{}-\eps}.
\end{align*}
Since $\eps>0$ was arbitrary we have 
\begin{equation*}
\inf_{F\in \hs(S), \left\| F\right\|{} =1}
\sup_{H\in \hs(S^*), \left\| H\right\|{}\le 1}
\left|\lma F,H \rma\right|
\le \frac{1}{\left\|\ler_S\right\|{} }
\end{equation*}
as required. 
\end{proof}

Thus $\left\|\ler_S\right\|{}$ measures efficiency of the pairing between $\hs(S)$ and 
$\hs(S^*)$.

If $S$ is the M\"obius image of a sphere then so is $S^*$, and using Remark \ref{R:anti-cr} we find that we are essentially pairing $\hs(S)$ with its conjugate; thus the pairing is perfectly efficient and $\left\|\ler_S\right\|{}=1$.

On the other hand, it follows from work of Bolt ([Bol1] , [Bol2]) that $\left\|\ler_S\right\|{}=1$ implies that $S$ is the M\"obius image of a sphere.

\begin{Remark} It follows from results in [BaLa] that 
\begin{equation*}
\left\|\ler_S\right\|{}\ge \max \left\{\left(1-|\belt_S(z)|^2\right)^{-1/2}\st z\in S\right\}
\end{equation*}
when $S$ is the smooth boundary of a strongly convex Reinhardt domain in $\C^2$. 
\end{Remark}

\end{document}